\documentclass[12pt]{article}
\usepackage{latexsym}
\usepackage{amsfonts}
\usepackage{amsmath}
\usepackage{amssymb}
\usepackage{graphicx}
\usepackage{latexsym}
\usepackage{amsfonts}
\usepackage{graphicx}
\usepackage{psfrag}

\input{tcilatex}
\begin{document}

\bigskip

\thispagestyle{empty}

\begin{center}
{\Large \textbf{On Lower bounds for variance and moments of unimodal
distributions}}

\vskip0.2inR. Sharma, R. Bhandari, R. Saini \\[0pt]
Department of Mathematics \& Statistics\\[0pt]
Himachal Pradesh University\\[0pt]
Shimla -5,\\[0pt]
India - 171005\\[0pt]
email: rajesh\underline{~~}hpu\underline{~~}math@yahoo.co.in
\end{center}

\vskip1.5in \noindent \textbf{Abstract. \quad }We provide an elementary
proof of the lower bound for the variance of continuous unimodal
distributions and obtain analogous bounds for the higher order central
moments. A lower bound for the $rth$ central moment of discrete distribution
is given and compared favorably with a related bound for discrete unimodal
distribution in literature.

\vskip0.5in \noindent \textbf{AMS classification \quad } 60E15.


\bigskip

\bigskip

\bigskip

\bigskip

\bigskip

\bigskip

\bigskip

\bigskip

\bigskip

\section{Introduction}

\setcounter{equation}{0} The arithmetic mean and variance of a continuous
random variable are respectively defined as%
\begin{equation*}
\mu _{1}^{^{\prime }}=\int\limits_{a}^{b}x\text{ }\phi \left( x\right) dx
\end{equation*}%
and%
\begin{equation*}
\sigma ^{2}=\int\limits_{a}^{b}\left( x-\mu _{1}^{^{\prime }}\right) ^{2}%
\text{ }\phi \left( x\right) dx
\end{equation*}%
where $\phi \left( x\right) $ is probability density function,%
\begin{equation*}
\int\limits_{a}^{b}\text{ }\phi \left( x\right) dx=1.
\end{equation*}%
A distribution is said to be unimodal at $x=M$ if $\phi \left( x\right) $ is
non-decreasing in $\left[ a,M\right) $ and non-increasing in $\left( M,b%
\right] .$ The special cases are non-increasing distributions $\left(
M=a\right) $ and non-decreasing distributions $\left( M=b\right) .$
Beginning with Gray and Odell $(1967),$ bounds for the variance have been
studied by several authors, see Jacobson $\left( 1969\right) $, Olshen and
Savage $\left( 1970\right) $, Seaman et al. $\left( 1987\right) $,
Dharmadhikari and Joag-Dev $\left( 1989\right) $ and Sharma and Bhandari (in
press).

An interesting inequality due to Johnson and Rogers $\left( 1951\right) $
says that for a unimodal distribution variance is bounded below by
(mean-mode)$^{2}/3$, that is 
\begin{equation*}
\sigma ^{2}\geq \frac{\left( \mu _{1}^{^{\prime }}-M\right) ^{2}}{3}.
\end{equation*}%
Jacobson $\left( 1969\right) $ gives a least upper bound for the variance of
unimodal distribution. The complementary inequality due to Jacobson $\left(
1969\right) $ says that%
\begin{equation*}
\sigma ^{2}\leq \frac{\left( b-a\right) ^{2}}{9}.
\end{equation*}%
Such simple and interesting inequalities are not popularly known. The proofs
of these inequalities are lengthy and tedious, see Dharmadhikari and
Joag-Dev $\left( 1989\right) $. The derivations given by Dharmadhikari and
Joag-Dev $\left( 1989\right) $ depends on the characterization of unimodal
distribution due to Shepp $\left( 1962\right) $ and, Olshen and Savage $%
\left( 1970\right) $. Recently, Sharma and Bhandari (in press) have given
elementary proofs of the upper bounds for the variance. In a similar spirit
we give here simple and elementary proof of the lower bound for the variance.

We also discuss lower bounds for the central moments of discrete
distribution. A discrete distribution $\{p_{1},p_{2},...,p_{n}\}$ with
support $\left\{ x_{1},x_{2},...,x_{n}\right\} $ such that $x_{1}\leq
x_{2}\leq ...\leq x_{n}$ is unimodal about $x=x_{k}$ if $p_{1}\leq p_{2}\leq
...\leq p_{k}\geq p_{k+1}\geq ...\geq p_{n},$ see Keilson and Gerber $\left(
1971\right) $ and Medgyessy $\left( 1972\right) $. A lower bound for the
discrete unimodal distribution $p_{i}$ with support $\left\{
...-2,-1,0,1,2,...\right\} $ and mode $M$ is%
\begin{equation}
3\sigma ^{2}\geq \left( \mu _{1}^{^{\prime }}-M\right) ^{2}+\left\vert \mu
_{1}^{^{\prime }}-M\right\vert ,  \label{1.1}
\end{equation}%
where%
\begin{equation*}
\mu _{1}^{^{\prime }}=\dsum\limits_{i=1}^{n}p_{i}x_{i}\text{\ and }\sigma
^{2}=\ \dsum\limits_{i=1}^{n}p_{i}\left( x_{i}-\mu _{1}^{^{\prime }}\right)
^{2},
\end{equation*}%
see Abouammoh and Mashhour $\left( 1994\right) $ and references therein.

We demonstrate here elementary proofs of the lower bounds for the variance
of continuous unimodal distributions ( Theorem 2.1-2.2, below). The
analogous bounds for the higher order moments are obtained (Theorem 2.3). We
prove a lower bound for the $r$th central moment of discrete distributions
and compare it favorably with the variance bounds given in Abouammoh and
Mashhour $\left( 1994\right) $ for discrete unimodal distributions.

\section{Main Results}

\setcounter{equation}{0} Our idea of the proof in the following theorems is
that the inequality of the type%
\begin{equation}
\int\limits_{a}^{b}\left( x-\alpha \right) \left( x-\beta \right) ~\phi
\left( x\right) dx\geq 0  \label{2.1}
\end{equation}%
yields a lower bound for the second order moment,%
\begin{equation}
\mu _{2}^{^{\prime }}\geq \left( \alpha +\beta \right) \mu _{1}^{^{\prime
}}-\alpha \beta .  \label{2.2}
\end{equation}%
For instance, the inequality \eqref{2.1} is always true for $\alpha =\beta ,$
from \eqref{2.2} we get that%
\begin{equation}
\mu _{2}^{^{\prime }}\geq 2\alpha \mu _{1}^{^{\prime }}-\alpha ^{2}.
\label{2.3}
\end{equation}%
The inequality \eqref{2.3} is valid for every real number $\alpha .$ The
function $f\left( x\right) =2cx-x^{2}$ attains its maximum at $x=c.$ The
inequality \eqref{2.3} therefore yields the classical lower bound $\mu
_{2}^{^{\prime }}\geq \mu _{1}^{^{\prime }2}$ at $\alpha =\mu _{1}^{^{\prime
}}.$ We show in the following theorem that when $\phi \left( x\right) $ is
non- increasing in $\left[ a,b\right] ,$ we can find distinct values of $%
\alpha $ and $\beta $ for which \eqref{2.1}\ holds and \eqref{2.2} gives a
better bound$.$

\vskip0.2in \textbf{Theorem 2.1. }Let $X$ be a random variable with
non-increasing distribution, $a\leq X\leq b$ $.$ Then%
\begin{equation}
\sigma ^{2}\geq \frac{\left( \mu _{1}^{^{\prime }}-a\right) ^{2}}{3}.
\label{2.4}
\end{equation}%
If $X$ has non-decreasing distribution,%
\begin{equation}
\sigma ^{2}\geq \frac{\left( \mu _{1}^{^{\prime }}-b\right) ^{2}}{3}.
\label{2.5}
\end{equation}%
\vskip0.2in\noindent\ \textbf{Proof.} We first consider the case when $\phi
\left( x\right) $ is non-increasing in $\left[ a,b\right] .$ For $a<\alpha
<\beta <b,$ we write%
\begin{equation}
\int\limits_{a}^{b}\left( x-\alpha \right) \left( x-\beta \right) ~\phi
\left( x\right) dx=\int\limits_{a}^{\beta }\left( x-\alpha \right) \left(
x-\beta \right) ~\phi \left( x\right) dx+\int\limits_{\beta }^{b}\left(
x-\alpha \right) \left( x-\beta \right) ~\phi \left( x\right) dx.
\label{2.6}
\end{equation}%
Since $\left( x-\alpha \right) \geq 0$ and $\left( x-\beta \right) \geq 0$
for $\beta \leq x\leq b,$%
\begin{equation}
\int\limits_{\beta }^{b}\left( x-\alpha \right) \left( x-\beta \right) ~\phi
\left( x\right) dx\geq 0.  \label{2.7}
\end{equation}%
Also $\left( x-\alpha \right) \leq 0$, $\left( x-\beta \right) \leq 0$ and $%
\phi \left( x\right) \geq \phi \left( \alpha \right) $ for $a\leq x\leq
\alpha $ and $\left( x-\alpha \right) \geq 0$, $\left( x-\beta \right) \leq
0 $ and $\phi \left( x\right) \leq \phi \left( \alpha \right) $ for $\alpha
\leq x\leq \beta ,$ therefore%
\begin{equation}
\int\limits_{a}^{\beta }\left( x-\alpha \right) \left( x-\beta \right) ~\phi
\left( x\right) dx\geq \phi \left( \alpha \right) \int\limits_{a}^{\beta
}\left( x-\alpha \right) \left( x-\beta \right) dx.  \label{2.8}
\end{equation}%
We conclude from \eqref{2.6}, \eqref{2.7} and \eqref{2.8} that%
\begin{equation}
\int\limits_{a}^{b}\left( x-\alpha \right) \left( x-\beta \right) ~\phi
\left( x\right) dx\geq \phi \left( \alpha \right) \int\limits_{a}^{\beta
}\left( x-\alpha \right) \left( x-\beta \right) dx.  \label{2.9}
\end{equation}%
It is clear from \eqref{2.9} that \eqref{2.1} holds if we choose $\alpha $
and $\beta $ such that%
\begin{equation}
\int\limits_{a}^{\beta }\left( x-\alpha \right) \left( x-\beta \right) dx=0.
\label{2.10}
\end{equation}%
We find from \eqref{2.10} that%
\begin{equation}
\alpha =\frac{2a+\beta }{3}.  \label{2.11}
\end{equation}%
Insert \eqref{2.11} into \eqref{2.2}, we see that%
\begin{equation}
\mu _{2}^{^{\prime }}\geq \frac{2}{3}\left( a+2\beta \right) \mu
_{1}^{^{\prime }}-\frac{2}{3}a\beta -\frac{\beta ^{2}}{3}.  \label{2.12}
\end{equation}%
We now find that value of $\beta $ for which \eqref{2.12} yields a greatest
lower bound. The function%
\begin{equation}
f\left( \beta \right) =\frac{2}{3}\left( a+2\beta \right) \mu _{1}^{^{\prime
}}-\frac{2}{3}a\beta -\frac{\beta ^{2}}{3}  \label{2.13}
\end{equation}%
with derivatives%
\begin{equation*}
f^{^{\prime }}\left( \beta \right) =\frac{4}{3}\mu _{1}^{^{\prime }}-\frac{2%
}{3}a-\frac{2\beta }{3}\text{ \ and \ }f^{\prime \prime }\left( \beta
\right) =-\frac{2}{3}
\end{equation*}%
has maximum at%
\begin{equation}
\beta =2\mu _{1}^{^{\prime }}-a.  \label{2.14}
\end{equation}%
Substitute value of $\beta $ from \eqref{2.14} in \eqref{2.12}, a simple
calculation leads to \eqref{2.4}, $\mu _{2}^{^{\prime }}=\sigma ^{2}+\mu
_{1}^{^{\prime }2}$. The inequality \eqref{2.5} follows on using similar
arguments. We have%
\begin{equation*}
\int\limits_{a}^{b}\left( x-\alpha \right) \left( x-\beta \right) ~\phi
\left( x\right) dx=\int\limits_{a}^{\alpha }\left( x-\alpha \right) \left(
x-\beta \right) ~\phi \left( x\right) dx+\int\limits_{\alpha }^{b}\left(
x-\alpha \right) \left( x-\beta \right) ~\phi \left( x\right) dx.
\end{equation*}%
If $\phi \left( x\right) $ is non-decreasing in $\left[ a,b\right] ,$%
\begin{equation*}
\int\limits_{a}^{\alpha }\left( x-\alpha \right) \left( x-\beta \right)
~\phi \left( x\right) dx\geq 0
\end{equation*}%
and%
\begin{equation*}
\int\limits_{\alpha }^{b}\left( x-\alpha \right) \left( x-\beta \right)
~\phi \left( x\right) dx\geq \phi \left( \beta \right) \int\limits_{\alpha
}^{b}\left( x-\alpha \right) \left( x-\beta \right) dx.
\end{equation*}%
Therefore%
\begin{equation*}
\int\limits_{a}^{b}\left( x-\alpha \right) \left( x-\beta \right) ~\phi
\left( x\right) dx\geq \phi \left( \beta \right) \int\limits_{\alpha
}^{b}\left( x-\alpha \right) \left( x-\beta \right) dx.
\end{equation*}%
The equation%
\begin{equation*}
\int\limits_{\alpha }^{b}\left( x-\alpha \right) \left( x-\beta \right) dx=0
\end{equation*}%
gives $3\beta =2b+\alpha $ and \eqref{2.2} becomes%
\begin{equation}
\mu _{2}^{^{\prime }}\geq \frac{2}{3}\left( b+2\alpha \right) \mu
_{1}^{^{\prime }}-\frac{2}{3}b\alpha -\frac{\alpha ^{2}}{3}.  \label{2.15}
\end{equation}%
The inequality \eqref{2.5} follows from \eqref{2.15}, $\alpha =b-2\mu
_{1}^{^{\prime }}.\blacksquare $

\vskip0.0in Note that for a non-increasing distribution $a\leq \mu
_{1}^{^{\prime }}\leq \frac{a+b}{2}.$ If $\mu _{1}^{^{\prime }}=a$, $\sigma
=0.$ If $\mu _{1}^{^{\prime }}=\frac{a+b}{2},$%
\begin{equation*}
\text{ }\phi \left( x\right) =\frac{1}{b-a}\text{\ and ~\ }\sigma ^{2}=\frac{%
\left( b-a\right) ^{2}}{12}.
\end{equation*}%
From \eqref{2.11} and \eqref{2.14}, $\alpha =\frac{2\mu _{1}^{^{\prime }}+a}{%
3}.$ Hence, for a non-increasing and non-uniform distribution $a<\alpha
<\beta <b.$The arguments given in the proof of Theorem 2.1 can now readily
be extended to prove the bounds for the unimodal distributions.

\vskip0.2in\noindent \textbf{\ Theorem 2.2.\ }Let $X$ be a random variable
such that $a\leq X\leq b$ and $X$ has a unimodal distribution at $M$, then%
\begin{equation}
\sigma ^{2}\geq \frac{\left( \mu _{1}^{^{\prime }}-M\right) ^{2}}{3}.
\label{2.16}
\end{equation}%
\vskip0.2in \noindent \textbf{Proof.} For the unimodal distribution, $\frac{%
a+M}{2}\leq \mu _{1}^{^{\prime }}\leq \frac{b+M}{2}.$ We first consider the
case when $M\leq \mu _{1}^{^{\prime }}\leq \frac{b+M}{2}.$ For $a\leq M\leq
\alpha \leq \beta \leq b,$ we write%
\begin{equation*}
\int\limits_{a}^{b}\left( x-\alpha \right) \left( x-\beta \right) ~\phi
\left( x\right) dx=\int\limits_{a}^{M}\left( x-\alpha \right) \left( x-\beta
\right) ~\phi \left( x\right) dx
\end{equation*}%
\begin{equation}
\text{ \ \ \ \ \ \ \ \ \ \ \ \ \ \ \ \ \ \ \ \ \ }+\int\limits_{M}^{\beta
}\left( x-\alpha \right) \left( x-\beta \right) ~\phi \left( x\right)
dx+\int\limits_{\beta }^{b}\left( x-\alpha \right) \left( x-\beta \right)
~\phi \left( x\right) dx.  \label{2.17}
\end{equation}%
It is easily seen on using arguments similar to those in the proof of
Theorem 2.1 that%
\begin{equation}
\int\limits_{a}^{M}\left( x-\alpha \right) \left( x-\beta \right) ~\phi
\left( x\right) dx\geq 0,\text{ }\int\limits_{\beta }^{b}\left( x-\alpha
\right) \left( x-\beta \right) ~\phi \left( x\right) dx\geq 0  \label{2.18}
\end{equation}%
and%
\begin{equation}
\int\limits_{M}^{\beta }\left( x-\alpha \right) \left( x-\beta \right) ~\phi
\left( x\right) dx\geq \phi \left( \alpha \right) \int\limits_{M}^{\beta
}\left( x-\alpha \right) \left( x-\beta \right) ~dx.  \label{2.19}
\end{equation}%
Combine \eqref{2.17}, \eqref{2.18} and \eqref{2.19}, we get%
\begin{equation*}
\int\limits_{a}^{b}\left( x-\alpha \right) \left( x-\beta \right) ~\phi
\left( x\right) dx\geq \phi \left( \alpha \right) \int\limits_{M}^{\beta
}\left( x-\alpha \right) \left( x-\beta \right) ~dx.
\end{equation*}%
Also,%
\begin{equation*}
\int\limits_{M}^{\beta }\left( x-\alpha \right) \left( x-\beta \right) ~dx=0
\end{equation*}%
gives%
\begin{equation}
\alpha =\frac{2M+\beta }{3}.  \label{2.20}
\end{equation}%
Therefore, \eqref{2.2} becomes%
\begin{equation*}
\mu _{2}^{^{\prime }}\geq \frac{2}{3}\left( M+2\beta \right) \mu
_{1}^{^{\prime }}-\frac{2}{3}M\beta -\frac{\beta ^{2}}{3}.
\end{equation*}%
The inequality \eqref{2.16} follows from the fact that the function%
\begin{equation*}
f\left( \beta \right) =\frac{2}{3}\left( M+2\beta \right) \mu _{1}^{^{\prime
}}-\frac{2}{3}M\beta -\frac{\beta ^{2}}{3}
\end{equation*}%
achieves its maximum at 
\begin{equation}
\beta =2\mu _{1}^{^{\prime }}-M  \label{2.21}
\end{equation}%
and $\sigma ^{2}=\mu _{2}^{^{\prime }}-\mu _{1}^{^{\prime }2}.$ For $M\leq
\mu _{1}^{^{\prime }}\leq \frac{b+M}{2}$ the values of $\alpha $ and $\beta $
obtained from \eqref{2.20} and \eqref{2.21} satisfy $a\leq M\leq \alpha \leq
\beta \leq b.$ On using similar arguments we find that \eqref{2.16} also
holds good when $\frac{a+M}{2}\leq \mu _{1}^{^{\prime }}\leq M$.$%
\blacksquare $

\vskip0.0in\noindent\ It is natural to consider the generalisation of
variance bounds for higher order central moments. We now show that the
inequality of the type%
\begin{equation*}
\int\limits_{a}^{b}f\left( x\right) \phi \left( x\right) dx\geq 0,\text{ \ }%
r=1,2,...,
\end{equation*}%
with%
\begin{equation}
f\left( x\right) =x^{r}-\frac{\beta ^{r}-\alpha ^{r}}{\beta -\alpha }x+\frac{%
\beta ^{r}\alpha -\alpha ^{r}\beta }{\beta -\alpha },\text{ }\alpha \neq
\beta  \label{2.22}
\end{equation}%
gives lower bound for $rth$ order moment,%
\begin{equation*}
\mu _{r}^{^{\prime }}\geq \frac{\beta ^{r}-\alpha ^{r}}{\beta -\alpha }\mu
_{1}^{^{\prime }}-\frac{\beta ^{r}\alpha -\alpha ^{r}\beta }{\beta -\alpha }
\end{equation*}%
where 
\begin{equation}
\mu _{r}^{^{\prime }}=\int\limits_{a}^{b}x^{r}\phi \left( x\right) dx.
\label{2.23}
\end{equation}%
For arbitrary distribution the limiting case $\left( \beta \rightarrow
\alpha \right) $ gives%
\begin{equation}
\mu _{r}^{^{\prime }}\geq r\alpha ^{r-1}\mu _{1}^{^{\prime }}-\left(
r-1\right) \alpha ^{r}.  \label{2.24}
\end{equation}%
The inequality \eqref{2.24} is valid for every real number $\alpha $ and
yields the greatest lower bound at $\alpha =\mu _{1}^{^{\prime }}.$ It
follows that $\mu _{r}^{^{\prime }}\geq \mu _{1}^{^{\prime }r}$ when $r$ is
even positive integer. Also, $\mu _{r}^{^{\prime }}\geq \mu _{1}^{^{\prime
}r}$ for every $r=1,2,...$ , $a\geq 0$.

\vskip0.0in\noindent\ An equivalent expression for polynomial \eqref{2.22} is%
\begin{eqnarray}
f\left( x\right) &=&x^{r}-\left( \alpha ^{r-1}+\alpha ^{r-2}\beta
+...+\alpha \beta ^{r-2}+\beta ^{r-1}\right) x  \notag \\
&&+\alpha \beta \left( \alpha ^{r-2}+\alpha ^{r-3}\beta +...+\alpha \beta
^{r-3}+\beta ^{r-2}\right) .  \label{2.25}
\end{eqnarray}%
Another equivalent form of our concern is given in the following lemma.

\vskip0.2in\noindent \textbf{\ Lemma 2.1. }For $x\geq 0,$ $\alpha \geq 0$
and $\beta \geq 0,$ the polynomial \eqref{2.22} can be written as 
\begin{equation}
f\left( x\right) =\left( x-\alpha \right) \left( x-\beta \right) g\left(
x\right)  \label{2.26}
\end{equation}%
where $g\left( x\right) \geq 0.$ If $r$ is even positive integer, %
\eqref{2.26} holds for all real $\alpha ,$ $\beta $ and $x$ with $g\left(
x\right) \geq 0$.

\vskip0.2in \noindent \textbf{Proof. }By Descarte's rule of sign if $r$ is
odd, $f\left( x\right) $ has only two positive real roots $\alpha $ and $%
\beta .$ The assertion of the lemma follows from the fact that $f\left(
x\right) >0$ for $0\leq x\leq \alpha $ and $x\geq \beta ,$ and $f\left(
x\right) <0$ for $\alpha \leq x\leq \beta .$ Note that $f\left( 0\right) =%
\frac{\beta ^{r}\alpha -\alpha ^{r}\beta }{\beta -\alpha }\geq 0.$ If $r$ is
even positive integer, \thinspace $\,f\left( x\right) $ has two real roots
and $g\left( x\right) $ is a polynomial of even degree whose roots occurs in
conjugate pairs. Hence $g\left( x\right) \geq 0$ for all real $%
x.\blacksquare $

\vskip0.2in\noindent \textbf{\ Theorem 2.3. }Let $X$ be a random variable
with non-increasing distribution, $a\leq X\leq b.$ Then%
\begin{equation}
\mu _{r}^{^{\prime }}\geq \frac{1}{2\left( r+1\right) }\frac{\left( 2\mu
_{1}^{^{\prime }}-a\right) ^{r+1}-a^{r+1}}{\mu _{1}^{^{\prime }}-a}.
\label{2.27}
\end{equation}%
\vskip0.2in \noindent \textbf{Proof. }It follows from Lemma 2.1 that%
\begin{equation}
\int\limits_{a}^{b}f\left( x\right) \phi \left( x\right)
dx=\int\limits_{a}^{b}\left( x-\alpha \right) \left( x-\beta \right) g\left(
x\right) \phi \left( x\right) dx,  \label{2.28}
\end{equation}%
where $g\left( x\right) \geq 0.$ It is evident that the integrand in %
\eqref{2.28} is positive or negative according as $\left( x-\alpha \right)
\left( x-\beta \right) $ is positive or negative. It follows on using
arguments similar to those used in the proof of Theorem 2.1 that%
\begin{equation}
\int\limits_{a}^{b}f\left( x\right) \phi \left( x\right) dx\geq \phi \left(
\alpha \right) \int\limits_{a}^{\beta }\left( x-\alpha \right) \left(
x-\beta \right) g\left( x\right) dx.  \label{2.29}
\end{equation}%
We choose $\alpha $ and $\beta $ such that%
\begin{equation}
\int\limits_{a}^{\beta }f\left( x\right) dx=0  \label{2.30}
\end{equation}%
From \eqref{2.22} and \eqref{2.30} we find that%
\begin{equation}
\alpha ^{r-1}+\alpha ^{r-2}\beta +...+\alpha \beta ^{r-2}=\frac{\left(
r-1\right) \beta ^{r}+\left( r-1\right) a\beta ^{r-1}-2a^{2}\beta
^{r-2}...-2a^{r-1}\beta -2a^{r}}{\left( r+1\right) \left( \beta -a\right) }.
\label{2.31}
\end{equation}%
It is easy to see from \eqref{2.29} - \eqref{2.31} that if $\alpha $ and $%
\beta $ satisfy \eqref{2.31}$,$%
\begin{equation}
\int\limits_{a}^{b}f\left( x\right) \phi \left( x\right) dx\geq 0.
\label{2.32}
\end{equation}%
Combine \eqref{2.25} and \eqref{2.32}, we get%
\begin{equation}
\mu _{r}^{^{\prime }}\geq \left( \alpha ^{r-1}+\alpha ^{r-2}\beta
+...+\alpha \beta ^{r-2}+\beta ^{r-1}\right) \mu _{1}^{^{\prime }}-\alpha
\beta \left( \alpha ^{r-2}+\alpha ^{r-3}\beta +...+\alpha \beta ^{r-3}+\beta
^{r-2}\right) .  \label{2.33}
\end{equation}%
Using \eqref{2.31} in \eqref{2.33},%
\begin{eqnarray}
\mu _{r}^{^{\prime }} &\geq &\left( \frac{\left( r-1\right) \beta
^{r}+\left( r-1\right) a\beta ^{r-1}-2a^{2}\beta ^{r-2}...-2a^{r-1}\beta
-2a^{r}}{\left( r+1\right) \left( \beta -a\right) }+\beta ^{r}\right) \mu
_{1}^{^{\prime }}  \notag \\
&&-\frac{\left( r-1\right) \beta ^{r}+\left( r-1\right) a\beta
^{r-1}-2a^{2}\beta ^{r-2}...-2a^{r-1}\beta -2a^{r}}{\left( r+1\right) \left(
\beta -a\right) }\beta .  \label{2.34}
\end{eqnarray}%
Let $h\left( \beta \right) $ denotes the right hand side expression of %
\eqref{2.34}. The function $h\left( \beta \right) $ with derivative%
\begin{equation*}
h^{^{\prime }}\left( \beta \right) =\frac{1}{r+1}\left( r\left( r-1\right)
\beta ^{r-2}+\left( r-1\right) \left( r-2\right) a\beta
^{r-3}+...+2a^{r-2}\right) \left( 2\mu _{1}^{^{\prime }}-\beta -a\right)
\end{equation*}%
has maximum at 
\begin{equation}
\beta =2\mu _{1}^{^{\prime }}-a.  \label{2.35}
\end{equation}%
It is lengthy to see that \eqref{2.34} implies \eqref{2.27} for $\beta =2\mu
_{1}^{^{\prime }}-a.$ Alternatively, substitute value of $\beta $ from %
\eqref{2.35} in \eqref{2.30}, we get%
\begin{equation}
\frac{\beta ^{r}-\alpha ^{r}}{\beta -\alpha }\mu _{1}^{^{\prime }}-\frac{%
\beta ^{r}\alpha -\alpha ^{r}\beta }{\beta -\alpha }=\frac{1}{2\left(
r+1\right) }\frac{\left( 2\mu _{1}^{^{\prime }}-a\right) ^{r+1}-a^{r+1}}{\mu
_{1}^{^{\prime }}-a}.  \label{2.36}
\end{equation}%
Also, from \eqref{2.32}, we have%
\begin{equation}
\mu _{r}^{^{\prime }}\geq \frac{\beta ^{r}-\alpha ^{r}}{\beta -\alpha }\mu
_{1}^{^{\prime }}-\frac{\beta ^{r}\alpha -\alpha ^{r}\beta }{\beta -\alpha }.
\label{2.37}
\end{equation}%
Combine \eqref{2.36} and \eqref{2.37}, we immediately get \eqref{2.27}.$%
\blacksquare $

\vskip0.2in\noindent \textbf{\ Theorem 2.4.\ }Let $X$ be a random variable
such that $a\leq X\leq b$ and $X$ has a unimodal distribution at $M$, then%
\begin{equation}
\mu _{r}^{^{\prime }}\geq \frac{1}{2\left( r+1\right) }\frac{\left( 2\mu
_{1}^{^{\prime }}-M\right) ^{r+1}-M^{r+1}}{\mu _{1}^{^{\prime }}-M}.
\label{2.38}
\end{equation}%
\vskip0.2in \noindent \textbf{Proof. }The proof follows easily\textbf{\ }on
using arguments similar to those used in the proofs of Theorem 2.2 -2.3.$%
\blacksquare $

\vskip0.0in It may be noted here that the bounds for the central moment%
\begin{equation*}
\mu _{2r}=\int\limits_{a}^{b}\left( x-\mu _{1}^{^{\prime }}\right) ^{2r}%
\text{ }\phi \left( x\right) dx
\end{equation*}%
follows from \eqref{2.38} on replacing $\mu _{1}^{^{\prime }}$ and $M$ by $0$
and $M-\mu _{1}^{^{\prime }},$ respectively$.$ For a unimodal distribution
at $x=M$ we have%
\begin{equation}
\mu _{2r}\geq \frac{1}{2\left( 2r+1\right) }\frac{\left( \mu _{1}^{^{\prime
}}-M\right) ^{2r+1}-\left( M-\mu _{1}^{^{\prime }}\right) ^{2r+1}}{\mu
_{1}^{^{\prime }}-M}=\frac{\left( \mu _{1}^{^{\prime }}-M\right) ^{2r}}{2r+1}%
.  \label{2.39}
\end{equation}%
For $r=1,$ the inequality \eqref{2.39} gives lower bound for the variance
proved in Theorem 2.2.

In case of a discrete unimodal distribution $p_{i}$ with prescribed support $%
\left\{ x_{1},x_{2},\ldots ,x_{n}\right\} $ and mean $\mu _{1}^{^{\prime }}$
we can find $\ j$ such that $x_{j-1}\leq \mu _{1}^{^{\prime }}\leq x_{j},$ $%
j=2,3,\ldots ,n.$ We give a lower bound for the central moment%
\begin{equation*}
\mu _{2r}=\dsum\limits_{i=1}^{n}p_{i}\left( x_{i}-\mu _{1}^{^{\prime
}}\right) ^{r}
\end{equation*}%
in terms of $x_{j-1},$ $x_{j}$ and $\mu _{1}^{^{\prime
}}=\dsum\limits_{i=1}^{n}p_{i}x_{i}.$ Also, see \cite{8}.

\vskip0.2in\noindent \textbf{\ Theorem 2.5. }Let $p_{i}$ be a discrete
distribution with support $\left\{ x_{1},x_{2},\ldots ,x_{n}\right\} $. If
the mean of the distribution is prescribed, $x_{j-1}\leq \mu _{1}^{^{\prime
}}\leq x_{j}$ and%
\begin{equation}
\mu _{2r}\geq \frac{\left( \mu _{1}^{^{\prime }}-x_{j-1}\right) \left(
x_{j}-\mu _{1}^{^{\prime }}\right) ^{2r}+\left( x_{j}-\mu _{1}^{^{\prime
}}\right) \left( \mu _{1}^{^{\prime }}-x_{j-1}\right) ^{2r}}{x_{j}-x_{j-1}},%
\text{ \ \ \ \ \ }x_{j-1}<x_{j}.  \label{2.40}
\end{equation}

\vskip0.2in \noindent \textbf{Proof. }For $x_{1}\leq \ldots \leq x_{j-1}\leq
x_{j}\leq \ldots \leq x_{n}$ all the $x_{i}$ lies outside $\left(
x_{j-1},x_{j}\right) ,j=2,\ldots ,n.$ It follows from Lemma 2.1 that%
\begin{equation}
x_{i}^{2r}-\frac{x_{j}^{2r}-x_{j-1}^{2r}}{x_{j}-x_{j-1}}x_{i}+\frac{%
x_{j}^{2r}x_{j-1}-x_{j-1}^{2r}x_{j}}{x_{j}-x_{j-1}}\geq 0.  \label{2.41}
\end{equation}%
The inequality \eqref{2.41} is true for any real number $x_{i}$ in $\left(
x_{j-1},x_{j}\right) .$Therefore, it must also hold for $x_{i}-\mu
_{1}^{^{\prime }}$ in $\left( x_{j-1}-\mu _{1}^{^{\prime }},x_{j}-\mu
_{1}^{^{\prime }}\right) $, that is%
\begin{eqnarray}
\left( x_{i}-\mu _{1}^{^{\prime }}\right) ^{2r} &\geq &\frac{\left(
x_{j}-\mu _{1}^{^{\prime }}\right) ^{2r}-\left( x_{j-1}-\mu _{1}^{^{\prime
}}\right) ^{2r}}{x_{j}-x_{j-1}}\left( x_{i}-\mu _{1}^{^{\prime }}\right) 
\notag \\
&&-\frac{\left( x_{j}-\mu _{1}^{^{\prime }}\right) ^{2r}\left( x_{j-1}-\mu
_{1}^{^{\prime }}\right) -\left( x_{j-1}-\mu _{1}^{^{\prime }}\right)
^{2r}\left( x_{j}-\mu _{1}^{^{\prime }}\right) }{x_{j}-x_{j-1}}.
\label{2.42}
\end{eqnarray}%
Also, the sum of the deviations of all the $x_{i}$ from the mean is zero, 
\begin{equation}
\dsum\limits_{i=1}^{n}p_{i}\left( x_{i}-\mu _{1}^{^{\prime }}\right) =0.
\label{2.43}
\end{equation}%
Multiplying both sides of \eqref{2.42} by $p_{i},$ add $n$ inequalities for $%
i=1,2,\ldots ,n$ and use \eqref{2.43}, we immediately get the inequality %
\eqref{2.40}.$\blacksquare $

\section{Numerical Examples}

\bigskip \setcounter{equation}{0} \textbf{1.}The distribution 
\begin{equation*}
\text{ }p_{i}=\left\{ 
\begin{array}{c}
-\frac{1}{4}i^{2}+\frac{1}{20}i+\frac{1}{2}\text{ \ \ for }i=-1,0,1 \\ 
0\text{ \ \ \ \ \ \ \ \ \ \ \ \ \ \ \ \ \ \ \ elsewhere}%
\end{array}%
\right. \text{ .}
\end{equation*}%
with support $x_{i}=\left\{ ...,-2,-1,0,1,2...\right\} $ is unimodal at $i=0.
$ The mean of the distribution is $\mu _{1}^{^{\prime }}=\frac{1}{10}.$ From %
\eqref{1.1} and \eqref{2.40}, we respectively have $\sigma ^{2}\geq \frac{11%
}{300}$ and $\sigma ^{2}\geq \frac{9}{100}.$ The inequality \eqref{2.40}
gives better estimates than \eqref{1.1}.

\vskip0.0in \textbf{2. }The distribution%
\begin{equation*}
\text{ }p_{i}=\left\{ 
\begin{array}{c}
\frac{4-\left\vert i\right\vert }{12}\text{ \ \ for }i=\pm 1,\pm 2,\pm 3 \\ 
0\text{ \ \ \ \ \ \ \ \ \ \ \ \ \ \ \ \ \ \ \ elsewhere}%
\end{array}%
\right. \text{ .}
\end{equation*}%
with support $x_{i}=\left\{ ...-15,-10,-5,5,10,15...\right\} $ is unimodal
at $i=-5.$ The mean of the distribution is $\mu _{1}^{^{\prime }}=0.$ The
inequality \eqref{1.1} is not applicable and classical bound $\mu
_{2r}^{^{\prime }}\geq \left( \mu _{1}^{^{\prime }}\right) ^{2r}$gives $\mu
_{2r}^{^{\prime }}\geq 0.$ From \eqref{2.40}, we have $\mu _{2r}^{^{\prime
}}\geq 5^{2r}.$

\vskip0.2in \textbf{Acknowledgements: }The support of the UGC-SAP is
acknowledged.


\begin{thebibliography}{99}
\bibitem{1} Abouammoh, A.M., Mashhour A.F. (1994). Variance upper bounds and
convolutions of $\alpha -$ unimodal distributions. \textit{Statist. Probab.
Lett.} 21:281-289.

\bibitem{2} Dharmadhikari, S.W., Joag-Dev, K. (1989). Upper bounds for the
variances of certain random variables.\textit{\ Comm. Statist. Theory
Methods }18\textit{:}3235-3247.

\bibitem{3} Gray, H.L., Odell, P.L. (1967). On least favorable density
functions. \textit{SIAM Review,}\textbf{\ }9\textbf{:}715-720.

\bibitem{4} Jacobson, H.I. (1969). The maximum variance of restricted
unimodal distributions.\textit{\ Ann. math. statist.,} 40:1746-1752.

\bibitem{5} Johnson, N. L., Rogers, C. A., (1951). The moment problem for
unimodal distributions. \textit{Ann. math. statist.} 22:433-439.

\bibitem{6} Keilson, J., Gerber, H. (1971). Some results for discrete
unimodality. \textit{J. Amer. Statist. Assoc}., 66:386-389.

\bibitem{7} Medgyessy, P. (1972). On unimodality of discrete distributions, 
\textit{Period. Math. Hungar.}, 2:245-257.

\bibitem{8} Muilwijk, J. (1966). Note on a theorem of M. N. Murthy and V. K.
Sethi. \textit{Sankhya Ser. B} 28:183.

\bibitem{9} Olshen, R.A., Savage, L.J., (1970). A generalized unimodality.%
\textit{\ J. Appl. Probab. }7:21-34.

\bibitem{10} Seaman, J.W., Young, D.M., Turner, D.W. (1987). On the variance
of certain bounded random variables.\textit{\ Math. Scientist,} 12:109-116.

\bibitem{11} Sharma, R., Bhandari, R. (2013). On variance upper bounds for
unimodal distributions, \textit{Comm. Statist. Theory Methods, }in press.

\bibitem{12} Shepp, L.A. (1962). Symmetric random walk.\textit{\ Trans.
Amer. Math. Soc.} 104:144-153.
\end{thebibliography}
\end{document}